\documentclass[
final,
nomarks
]{dmtcs-episciences}

\usepackage[utf8]{inputenc}
\usepackage{subfigure}
\usepackage{ amssymb }

\usepackage{graphicx}

\usepackage[round]{natbib}

\author{Samvel Kh. Darbinyan}
\title [On pre-Hamiltonian and shorts cycles in Digraphs] {Sufficient conditions for a digraph to contain:
a pre-Hamiltonian cycle and cycles of lengths
3 and 4}
\affiliation{Institute for Informatics and Automation Problems of NAS RA \\Yerevan, Armenia}
\keywords{Digraph, Cycle, Hamiltonian cycle, pre-Hamiltonian cycle, semi-degree}

\begin{document}
\publicationdata
{vol. 28:3}
{2026}
{12}
{10.46298/dmtcs.16867}
{2025-11-06; 2025-11-06; 2026-06-05}
{2026-09-12}

\maketitle

\begin{abstract}
 \vspace{0.2cm} 

{Let $D$ be a digraph of order $p\geq5$ with minimum degree at least $p-1$ and with minimum semi-degree at least $p/2-1$.
 In his excellent and renowned paper, ``Long Cycles in Digraphs" (Proc. London Mathematical Society (3), 42 (1981)), Thomassen fully characterized the following for $p=2n+1$: (i) $D$ has a cycle of length at least $2n$; and (ii) $D$ is Hamiltonian. Motivated by this result, and building on some of the ideas in Thomassen's paper,  we investigated the Hamiltonicity  (when $p$ is even) and pancyclicity (when $p$ is arbitrary) of such digraphs.  We  have given a complete description of whether such digraphs are Hamiltonian ($p$ is even), and  pancyclic ($p$ is arbitrary). Since the proof is very long, we have divided it into three parts.
 In this paper, we  provide a full description of the following: (iii)  for $k=3$ and $k=4$,  the digraph $D$ contains a cycle of length $k$; and (iv) the digraph $D$ contains a pre-Hamiltonian cycle, i.e., a cycle of length $p-1$.} 
\end{abstract}

\section{Introduction}

In this paper, we consider finite digraphs without loops and multiple arcs, which may contain opposite arcs with the same end-vertices, i.e., a cycle of length two. We shall assume that the reader is familiar with the standard terminology on digraphs and graphs (undirected graphs). 
We refer the reader to  \cite{[1]}  for digraphs, and \cite{[27]}) for graphs, for the terminology and notation not defined in this paper.
All cycles and paths are assumed to be simple and directed. A cycle (respectively, path) in  a digraph $D$ is called a {\em Hamiltonian cycle} (respectively, a {\em Hamiltonian path}) if it contains  all the vertices of $D$  and, if $D$ has  a Hamiltonian cycle, we say that $D$ is {\em Hamiltonian}. 
 A cycle in a digraph $D$ of order $p\geq 3$  is called a {\em pre-Hamiltonian cycle}  if it contains exactly $p-1$ vertices of $D$. 
Further digraph terminology and notation are given in the next section. 
 One of the most fundamental and widely studied problems in digraph theory is finding sufficient conditions for a digraph to contain a Hamiltonian cycle or cycles with previously given properties. 
 
 There are many sufficient conditions for the existence of a Hamiltonian cycle in  digraphs (see, e.g., \cite{[2]}, \cite{[3]}, \cite{[20]}, \cite{[9],[10]}, \cite{[11]},  \cite{[17]}, \cite{[22]}, \cite{[23]}, \cite{[24]} and \cite{[28]}). For more information on Hamiltonian digraphs, see the survey papers \cite{[4]} and \cite{[21]}, as well as the book by \cite{[1]}. 

 \cite{[17]} proved that every strong digraph on $p$ vertices and with minimum degree at least $p$ is Hamiltonian. 
 \cite{[24]} observed that  a digraph of order $p$ with minimum semi-degree at  least $p/2$ is Hamiltonian and  raised  the problem of describing all the strong non-Hamiltonian  digraphs of order $p$ and minimum degree  $p-1$. In his excellent and renowned paper, \cite{[26]}  provided a structural characterization of the Nash-Williams problem. 

 Let $D$ be a digraph  of order  $p\geq 5$  with  minimum degree at least $p-1$ and with minimum semi-degree  at least $p/2-1$. In the same paper, \cite{[26]}  fully characterized the following for $p=2n+1$: (i) $D$ has a cycle of length at least $2n$; and (ii) $D$ is Hamiltonian.
 
Motivated by the last results, we investigated the Hamiltonicity  (when $p$ is even) and pancyclcity (when $p$ is arbitrary) of such digraphs. Building on some of the ideas in Thomassen's paper, we  have given a complete description of whether such digraphs are Hamiltonian ($p$ is even) and  pancyclic ($p$ is arbitrary). Since the proof is very long, we have divided it into three parts. Note that these results were only reported in \cite{[6]} and in \cite{[7]}, without their proofs.  The proof have never been published.  In this paper, we provide a full description of the following: (i)  for $k=3$ and $k=4$,  the digraph $D$ contains a cycle of length $k$; and (ii) the digraph $D$ contains a pre-Hamiltonian cycle, i.e., a cycle of length $p-1$.

 \section {Terminology and notation}
 
We denote the vertex set and arc set of a digraph $D$ by $V(D)$ and $A(D)$, respectively.  The {\em order} of a digraph $D$ is the number of its vertices. The arc of a digraph $D$ directed from $x$ to $y$ is denoted by  $xy$ or $x\rightarrow y$.
   If $xy$ is an arc, then we say that $x$ {\em dominates} $y$ (or $y$ is {\em dominated} by $x$). For  a pair of subsets $F$ and  $B$ of $V(D)$, we define $A(F\rightarrow B)$  as the set $\{xy\in A(D) | x\in F, y\in B\}$ and $A(F,B)=A(F\rightarrow B)\cup A(B\rightarrow F)$. If  $x\in V(D)$ and $F=\{x\}$, we often write $x$ instead of $\{x\}$. For disjoint subsets  $F$ and $B$  of $V(D)$,  $F\rightarrow B$ means that every vertex of $F$ dominates every vertex of $B$. If $S\subset V(D)$, $F\rightarrow B$ and $B\rightarrow S$, then we write $F\rightarrow B\rightarrow S$ (for short). The {\em out-neighborhood} of vertex $x$ is the set $O(x)=\{y\in V(D) | xy\in A(D)\}$ and $I(x)=\{y\in V(D) | yx\in A(D)\}$ is the {\em in-neighborhood} of $x$. Similarly, if $F\subseteq V(D)$ then $O(x,F)=\{y\in F | xy\in A(D)\}$ and $I(x,F)=\{y\in F | yx\in A(D)\}$. The {\em out-degree} of $x$ is $od(x)=|O(x)|$ and $id(x)=|I(x)|$ is the {\em in-degree} of $x$. We call the out-degree and in-degree of a vertex its {\em semi-degrees}.
Similarly, $od(x,F)=|O(x,F)|$ and $id(x,F)=|I(x,F)|$. The {\em degree} of the vertex $x$ in $D$ is defined as $d(x)=id(x)+od(x)$. A digraph $D$ is $k$-regular if, for any vertex $x\in V(D)$, $d(x)=k$. The subdigraph of $D$ induced by a subset $F$ of $V(D)$ is denoted by $D[F]$. 
     A {\em dipath} (for short, {\em path})  is a digraph with the  vertex set $\{x_1,x_2,\ldots ,x_m\}$ ($m\geq 2$) and the arc set  $\{x_ix_{i+1}\, |\, i\in [1,m-1]\}$. This path is called   a path from $x_1$ to $x_m$ or is an $(x_1,x_m)$-path and is denoted by $x_1x_2\ldots x_m$. The {\em dicycle} (for short, {\em  cycle})  $x_1x_2\ldots x_nx_1$ is the digraph obtained from the path $x_1x_2\ldots x_n$ by adding the arc $x_nx_1$.
  The {\em length} of a cycle or a path is the number of its arcs. A {\em $k$-cycle} is a cycle of length $k$.  $C_k$  denotes a cycle of length $k$. 
 For a cycle  $C_k:=x_1x_2\ldots x_kx_1$, the indices are taken modulo $k$, i.e., $x_s=x_i$ for every $s$ and $i$ such that  $i\equiv s\, \hbox{mod} \,k$.
  
 If $B\subset V(D)$ or $B$ is a subdigraph of a digraph $D$, then we write $D-B$ for  $D[V(D)-B]$ and $D[V(D)-V(B)]$.   Let $x,y\in V(D)$ be two distinct vertices. We write $a^+(x,y)=1$, if $xy\in A(D)$ and  $a^+(x,y)=0$  otherwise, and $a(x,y)$ means  the  number of arcs between the vertices $x$ and $y$, in particular, $a(x,y)=0$ means that the vertices $x$ and $y$ are not adjacent.
For integers $a$ and $b$, $a\leq b$, by $[a,b]$ we denote the set $\{a,a+1,\ldots , b\}$. 

A digraph $D$ is {\em strong}  if for every pair  of distinct vertices $x$, $y$ of $D$ there exists an $(x,y)$-path and a $(y,x)$-path. 
  A digraph $D$ is {\em $k$-strong}, $k\geq 1$, if $|V(D)|\geq k+1$ and $D-A$ is strong for any set $A$ of at most $k-1$ vertices. 
 The largest integer $k$ such that $D$ is $k$-strong is the {\em vertex-strong connectivity number} of $D$ (denoted by $k(D))$.
Two distinct vertices $x$ and $y$ in $D$ are {\em adjacent} if $xy\in A(D)$ or $yx\in A(D) $ (or both).
   If $xy$ is an arc of a digraph $D$, then we say that $x$ {\em dominates} $y$. An empty digraph is a digraph that contains no arcs. A set $S$ of arcs in a digraph $D$ is called independent if no two arcs of $S$ have a common vertex.
 For a graph (undirected graph) $G$, we denote by $G^*$ the symmetric digraph obtained from $G$ by replacing every edge $xy$ with the pair $xy, yx$ of arcs.  
 The undirected cycle   of length $k$ in $G$, is denoted by $Q_k$. (As usual, we use $C_k$ to represent a directed cycle of length $k$). 
 A digraph $D$ (undirected graph $G$) of order $p \geq 3$ is pancyclic if it contains a directed (undirected) cycle of length $k$ for every $k\in [3, p]$. 
 $K_n$ denotes the complete undirected graph with $n$ vertices, and $K_{p,q}$ denotes the complete undirected     bipartite graph with partite sets of cardinalities $p$ and $q$. 
 If $G_1$ and $G_2$ are undirected graphs, then $G_1\cup G_2$ is the disjoint union of $G_1$ and $G_2$ and $G_1+G_2$ is the join of $G_1$ and $G_2$. 
 The converse digraph of a digraph $D$ is the digraph obtained from $D$ by reversing the direction of all arcs, and is denoted by $D^{rev}$. 
We will use  {\em the principle of digraph duality}:  A digraph $D$ contains a subdigraph $H$ if and only if the converse digraph of $D$ contains the converse digraph of  $H$.

\section{ Preliminaries and additional notation}

Let us recall some well-known lemmas used in this paper.\\

\noindent\textbf{Lemma 3.1.} 
  (\cite{[18]}). Let $D$ be a digraph of order $p\geq 3$  containing a cycle $C_m$ of length $m$ with $m\in [2,p-1]$. Let $x$ be a vertex not contained in this cycle. If $d(x,V(C_m))\geq m+1$, then for every $k\in [2,m+1]$, $D$ contains a cycle $C_k$including $x$.\\

The following lemma, which  is a slight modification of a lemma by  \cite{[5]},  will be used extensively in the proofs our results.

\noindent\textbf{Lemma 3.2.}  Let $D$ be a digraph of order $p\geq 3$  containing a path $P:=x_1x_2\ldots x_m$ with $m\in [2,p-1]$ and let $x$ be a vertex not contained in this path. If one of the following holds:

 (i) $d(x,V(P))\geq m+2$; 

 (ii) $d(x,V(P))\geq m+1$ and $xx_1\notin A(D)$ or $x_mx_1\notin A(D)$; 

 (iii) $d(x,V(P))\geq m$, $xx_1\notin A(D)$ and $x_mx\notin A(D)$;
\noindent\textbf{}then there is an  $i\in [1,m-1]$ such that $x_ix,xx_{i+1}\in A(D)$, i.e., $D$ contains a path $x_1x_2\ldots x_ixx_{i+1}\ldots x_m$ of length $m$  (we say that  $x$ can be inserted into $P$ or the path  $x_1x_2\ldots x_ixx_{i+1}\ldots x_m$ is an extended path obtained from $P$ with  $x$).\\ 

The following is an immediate consequence of Lemma 3.2:

\noindent\textbf{Lemma 3.3.}  Let $D$ be a digraph of order $p\geq 4$  and let  $P:=x_1x_2\ldots x_m$ with $m\in [2,p-2]$ be a path of  maximal length from $x_1$ to $x_m$ in $D$. If the induced subdigraph $D[V(D)\setminus V(P)]$ is strong and  $d(x,V(P))=m+1$ for every vertex $x\in V(D)\setminus V(P)$, then there is an integer $l\in [1, m]$ such that $O(x,V(P))=\{x_1,x_2,\ldots ,x_l\}$ and $I(x,V(P))=\{x_l,x_{l+1},\ldots , x_m\}$.\\

 We will now introduce the following notation:\\

 \noindent\textbf {Notation 1.} {\em For every integer $n\geq 1$, let $\cal{H}$$(n,n)$ denote the set of digraphs $D$ of order $2n$  such that 
$V(D)=F\cup B$,\, $D[F] \cong D[B]\cong K_{n}^*$,  there is no arc from a vertex of $B$ to a vertex of $F$,  and for every vertex $x\in F$ (respectively, $y\in B$) $d^+(x,B)\not=0$ (respectively, $d^-( y,F)\not=0$)}.\\

\noindent\textbf {Notation 2.} {\em For every integer $n\geq 2$, let $\cal{H}$$(n,n-1,1)$ denote the set of digraphs $D$ of order $2n$ such that  $V(D)= F\cup B\cup \{a\}$, $|F|=|B|+1=n$,\, $D[ F]$ contains no arc, $D$ contains all the possible arcs between $F$ and $B$, 
   $D[B\cup \{a\}]$ is an arbitrary digraph of order  $n$  such that in $D$ either $I(a)=B$ and   $a\rightarrow F$ or $O(a)=B$ and   $F\rightarrow a$}.\\
   
 \noindent\textbf {Notation 3.} {\em For every integer $n\geq 2$ define a digraph $H_{2n}$ of order $2n$ as follows: \,$V(H_{2n})=F\cup B\cup \{x,y\}$, \, $H_{2n}[F\cup \{x\}] \cong H_{2n}[B\cup \{y\}] \cong K_{n}^*$, there is no arc between $F$ and $B$,  $O(x)=F\cup \{y\}$, $I(y)=B\cup \{x\}$ 
 and $I(x)=O(y)=F\cup B$}. 
$H'_{2n}$ is a digraph obtained from $H_{2n}$ by adding the arc $yx$ (see Figure 1(a)).
\begin{figure*}[!h]
\begin{center}
\includegraphics[width=14cm]{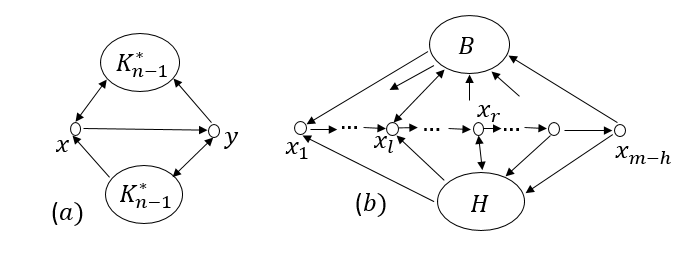}
\vspace{-5mm}
 \caption{(a) The digraph $H_{2n}$, (b) an illustration  for (8) in the proof of Theorem 5.1.}
\label{fig1}
\end{center}
\end{figure*}

It is not difficult to check that if
  $D\in \cal{H}$$(n,n)$ $\cup \, \cal{H}$$(n,n-1,1)\cup \{H_{2n},H'_{2n}\}$, then $D$ is not Hamiltonian and its each vertex   has degree   at least $2n-1$ and its semi-degrees  are at least $n-1$. Moreover, also $D^{rev}\in \cal{H}$$(n,n)$ $\cup \, \cal{H}$$(n,n-1,1)\cup \{H_{2n},H'_{2n}\}$.\\

\noindent\textbf{Notation 4.} By $\cal{B}$${(n,n)}$ we denote the set of  balanced bipartite digraphs of order $2n\geq 6$, obtained from $K^*_{n,n}$ by deleting $k$ independent arcs, where $0\leq k\leq n$. \\

Note that if $D\in$\,$\cal{B}$${(n,n)}$, then it is Hamiltonian, contains a cycle of length 4 and contains no cycle of lengths 3 and $2n-1$.\\

 \noindent\textbf{Remark 1.} Let $D$ be a digraph of order $2n\geq 6$ with minimum degree at least $2n-1$ and with minimum semi-degree at least $n-1$. If $V(D)=F\cup H$ such that $|F|=|H|$ and $D[F]$, $D[H]$ are empty digraphs, then $D\in$$\,\,\cal{B}$${(n,n)}$.\\

It is easy to check that the following lemma is true.

\noindent\textbf{Lemma 3.4.} Let $D$ be a digraph of order $p\geq 3$    with  minimum degree  at least $p-1$ and with minimum in- and out-degrees  at least $p/2-1$. Then

(i)  either $D$ is strong or $p=2n$ and $D\in \cal{H}$$(n,n)$;

(ii) if $B\subset V(D)$,  $|B|\geq (p+1)/2$ and  $x\in V(D)\setminus B$, then  $A(x\rightarrow B)\not=\emptyset$ and $A(B\rightarrow x)\not=\emptyset$. 
Moreover, if $|B|=\lfloor p/2\rfloor$ and  $A(x\rightarrow B)=\emptyset$  (respectively, $A(B\rightarrow x)=\emptyset$),  then $x\rightarrow V(D)\setminus (B\cup \{x\})$ (respectively, $V(D)\setminus (B\cup \{x\}\rightarrow x$). 

 \section {A sufficient condition for the existence of    cycles of lengths  3 and 4 in a digraph.}

The next two results will be used in the proof of Theorem 4.3.\\

\noindent\textbf{Theorem 4.1.} (\cite{[19]}). Let $G$ be an undirected graph of order $2n+1\geq 7$  with minimum degree  at least $n$. Then precisely one of the following holds: (a) $G$ is pancyclic;    (b) $ G\cong (K_n\cup K_n)+K_1$; or (c) $ K_{n,n+1}\subseteq G \subseteq K_n+\overline K_{n+1}$, where $\overline K_{n+1}$ is  the empty digraph with $n+1$ vertices.\\

\noindent\textbf{Theorem 4.2.}  (\cite{[25]}). Let $D$ be a strong digraph of order $p\geq 3$. If for every pair $x,y$ of nonadjacent  distinct vertices $d(x)+d(y)\geq 2p$, then $D$ is pancyclic or $p$ is even and $D\cong K_{p/2,p/2}^*$.\\

 Now we define the digraphs  $H^{'}_{6}$,  $H^{''}_{6}$ and  $H^{'''}_{6}$ as follows: \\

 (i) Let $Q_6$ be an undirected cycle $xvwuyzx$ of length 6. We denote by $H'_6$  the digraph obtained from $Q_6^*$ by adding the arcs $xy$, $xw$, $zu$ and $vu$.
  (See Figure 2 (a)).\\ 
 
(ii) Let $C_6=xyvzwux$ be a cycle of length 6. We denote by $H''_6$  the digraph obtained from $C_6$ by adding the arcs $xv$, $vx$, $xu$, $xw$, $uz$, $yu$, $yz$, $zy$, $wz$ and $wv$ (see Figure 3 (a)).\\

(iii) Let $C_6=xwvzyux$ be a cycle of length 6. We denote by $H'''_6$ the digraph obtained from $C_6$ by adding the arcs $xy$, $yv$, $vx$, $uz$, $zw$, $wu$, $xu$, $yz$ and  $vw$.  (See Figure 2 (b)). In all figures, an undirected edge represents two  arcs of opposite directions.\\

\begin{figure*}[!h]
\begin{center}
\includegraphics[width=11cm]{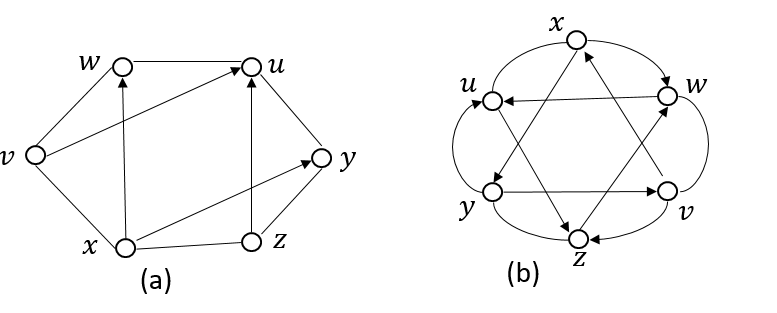}
\vspace{-5mm}
\caption{(a) The digraph $H'_6$ and  (b) the digraph $H'''_6$.}
\label{figure 2}
\end{center}
\end{figure*}


\noindent\textbf{Theorem 4.3.} Let $D$ be a digraph of order $p\geq 5$   with  minimum degree  at least $p-1$ and with  minimum semi-degree   at least $p/2-1$. Then the following holds: 

(i) either $D$ has a cycle of length 3 
or $D\in \cal{B}$$(n,n)$\,\,$\cup \,\{Q_5^*,  K_{n,n+1}^* \}$, where $n=\lfloor p/2 \rfloor$. 

(ii) either $D$ has a cycle of length 4  or $D\in \cal{H}$$(3,3)$ $\cup \{ Q_5^*, H^{'}_6, H^{''}_6, H^{'''}_{6},  [(K_2 \cup K_2)+K_1]^*\}$.

\begin{proof}  Suppose that  $D$ is a symmetric digraph, i.e.,   from $xy\in A(D)$ it follows that $yx\in A(D)$. Therefore,  $id(x)=od(x)$ for all vertices of $D$.  Then  $id(x)=od(x)\geq n$ since $d(x)\geq p-1$. If $p=2n$, then $D$ satisfies the conditions of Theorem 4.2. Therefore, either $D$ is pancyclic or $D= K^*_{n,n}\in \cal{B}$$(n,n)$. If $p=2n+1\geq 7$, then the underlying graph $G$ of $D$ satisfies the conditions of Theorem 4.1. 
Therefore, either (a) $G$ is pancyclic or    (b) $ G\cong (K_n\cup K_n)+K_1$ or (c) $ K_{n,n+1}\subseteq G \subseteq K_n+\overline K_{n+1}$. If (a) or (b) holds,  then it is clear that $D$ contains  cycles of lengths 3 and 4.
In   case (c),  $D$  contains a cycle of length 4, and contains no cycle of length 3 if and only if $D=K^*_{n,n+1}$. 
Thus we have, if $D$ is a symmetric digraph of order at least 6, then the theorem is true. Let now  $p=2n+1 =5$. Then it is not difficult to show that $D$ contains no cycle of length 3 if and only if $D\in \{Q^*_5,K^*_{2,3})\}$ 
and $D$ contains no cycle of length four if and only if $D\in \{Q^*_5,[(K_2 \cup K_2)+K_1]^*\}$.

In the following, suppose that $D$ is not symmetric. Then   $D\notin \{ Q_5^*, K_{n,n+1}^*, K^*_{n,n}, [(K_2 \cup K_2)+K_1]^*\}$. \\

(i)  Assume, to the contrary,  that $D$ has no cycle of length 3. Then
 for any arc $xy\in A(D)$ we have
$$
O(y)\cap I(x)=\emptyset. \eqno (1)
$$
Let $uv$ be an arbitrary arc in $D$ such that $vu\notin A(D)$. Then from (1),
$od(v)\geq p/2-1$ and $id(u)\geq p/2-1$ it follows that $p\geq |O(v)|+|I(u)|+2\geq p/2-1+p/2-1+2=p$.
Therefore, $p=2n\geq 6$ and $od(v) =id(u)= n-1$. Hence we have $id(v)\geq n$ and $od(u)\geq n$ since $d(x)\geq 2n-1$ for every vertex $x$ of $D$. Since $uv$ is an arbitrary arc in $D$, we can assume that $u\rightarrow I(u)$ and $O(v)\rightarrow v$. Therefore, $D[ O(v)]$ and $D[I(u)]$ are empty digraphs.
This together with $n\geq 3$ and (1) implies that for each vertex $z_1\in O(v)$ and $z_2\in I(u)$  we have $A(z_1\rightarrow I(u))\not=\emptyset$
and $A(O(v)\rightarrow z_2)\not=\emptyset$.  Therefore, $A(u\rightarrow O(v))=A(I(u)\rightarrow v)=\emptyset$ since $D$ contains no cycle of length 3. Thus we have $d(u, O(v))=d(v,I(u))=0$,  and hence the digraphs 
$D[\{u\}\cup O(v)])$ and $D[\{v\}\cup I(u)])$ are empty digraphs.
Then by Remark 1,  $D\in$ $\cal{B}$ ${(n,n)}$.\\

(ii) Assume   that $D$ contains no cycles of length 4, For each arc $xy\in A(D)$ put 
$$
S(x,y):=I(x)\cap O(y) \quad \hbox {and} \quad  E(x,y):=V(D)\setminus (O(y)\cup I(x)\cup \{x,y\}).
$$

Since $D$ has no  cycle of length 4, we see that
$$
A(O(y)\setminus \{x\} \rightarrow I(x)\setminus \{y\})=\emptyset. \eqno (2)
$$

 Let us consider the following cases.

\noindent\textbf {Case 1.} There exists an arc $xy\in A(D)$  such that $yx\notin A(D)$ and $od(y)\geq n$ or $id(x)\geq n$.

By the digraph duality,  we can assume that $od(y)\geq n$.  Then from (2) and Lemma 4(ii) it follows that 
$$
I(x)\subseteq O(y),\,\, \hbox {i.e.}, \,\, I(x)=S(x,y) \, \,\hbox {and} \quad A(D[S(x,y)])=\emptyset. \eqno (3)
$$

First, let us  prove the following claim.

\noindent\textbf{Claim 4.1.} $I(x)= O(y)$. 

\begin{proof} (of Claim 4.1) Assume that Claim 4.1 is not true.
 Then by (3), $O(y)\setminus I(x)\not= \emptyset$. Let $z$ be an arbitrary vertex in $O(y)\setminus I(x)$. By (2), we have $A(z\rightarrow \{x\}\cup I(x))=\emptyset$. From this  and Lemma 3.4(ii) it follows that  $p/2\leq |\{x\}\cup I(x)| \leq \lfloor p/2 \rfloor$, $p=2n\geq 6$, $id(x)=n-1$ and  
$$
O(z)=V(D)\setminus (\{x,z\}\cup I(x)),   \eqno {(4)}
$$
in particular, $zy\in A(D)$, $od(z)=n-1$,  $id(z)\geq n$ and $z\rightarrow E(x,y)$.     If $xz\in A(D)$, then $C_4=xzyux$, where $u\in S(x,y)$, which contradicts    our assumption. Therefore $xz\notin A(D)$. Then from $id(x)=n-1$, $id(z)\geq n$ and Lemma 3.4(ii) it follows that there exists $u\in I(x)$ such that  $uz\in A(D)$. Then by (4), it is not difficult to see that $D[O(y)\setminus I(x)]$ is an empty digraph,  and hence $O(y)\setminus I(x)=\{z \}$, $|E(x,y)|=n-2\geq 1$. Notice that $A(E(x,y)\rightarrow y)=\emptyset$, for otherwise, $C_4=yuzwy$, where $wy\in A(E(x,y)\rightarrow y)$.\\

Now, we will prove Claim 4.1 by examining two cases (a) and (b).\\

\noindent\textbf { Case (a).} For each vertex $v\in I(x)$ there is a vertex $v'\in E(x,y)$ such that $v'v\in A(D)$.

 If there exists  a vertex $v\in I(x)$ such that $vy\in A(D)$, then $C_4=vyzv'v$, a contradiction. We may therefore assume that  $A(I(x)\rightarrow y)=\emptyset$.
 Thus we have, $A(I(x)\cup E(x,y)\rightarrow y)=\emptyset$.  Then by Lemma 3.4(ii), $|I(x)|+|E(x,y)|=n-1+n-2\leq n$, which in turn implies that  $n=3$ and $|E(x,y)|=n-2= 1$.    Let $E(x,y):=\{w\}$ and $I(x):=\{u,v\}$. Under the condition of Case (a), we have $w \rightarrow \{u,v\}$. Note that $a(u,v)=a(x,z)=0$. 
 
 Assume first that $wz\in A(D)$. (See Figure 3 (b)) for an illustration.) Then it is easy to see that $A(\{u,v\}\rightarrow w)=\emptyset$.  From $a(u,v)=0$ and  $A(z\rightarrow \{u,v\})=A(\{u,v\}\rightarrow \{w,y\})=\emptyset$ it follows that $\{u,v\}\rightarrow z$, $od(u)=od(v)=2$, $id(u)=id(v)=3$ and $x\rightarrow \{u,v,w\}$. 
 Therefore, $D \cong H^{''}_{6}$.
 \begin{figure*}[!h]
\begin{center}
\includegraphics[width=12cm]{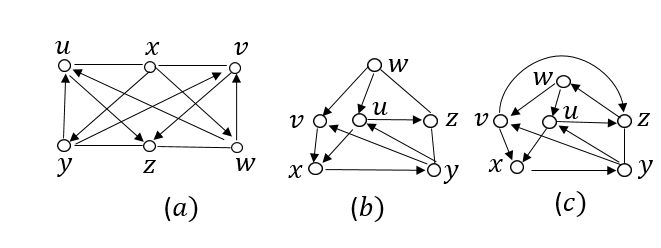}
\vspace{0 cm}
 \caption{(a) The digraph $H''_6$,  and  (b),  (c) an illustrations for  Case (a)} {in the proof of Theorem 4.3.}
\label{fig4}
\end{center}
\end{figure*}

Assume second that $wz\notin A(D)$. Since $A(w \rightarrow \{x,y,z\})=\emptyset$, 
we have $d^+(w)=2$ and $d^-(w)\geq 3$.  From $A(\{x,w\}\rightarrow z)=A(z\rightarrow \{u,v,x\})=\emptyset$ it follows that $I(z)=\{u,v,y\}$ and $O(z)=\{w,y\}$. (See  Figure 3 (c) for an illustration.) Assume that $xu\in A(D)$. Then $A(u\rightarrow  \{y,v,w\})=A(\{v,z\}\rightarrow  u)=\emptyset$. Since $A(\{y,u\}\rightarrow w)=\emptyset$ and $d^-(w)\geq 3$ it follows that $\{x,v\}\rightarrow  w$ and $xv\notin A(D)$. Therefore, $C_6=zyuxwvz$ is a cycle in $D$, the arcs $yv$, $vx$, $xy$, $uz$, $zw$, $wu$, $vw$, $xu$, $yz$ are in $A(D)$ and $D$ contains no other arcs. This means that  $D\cong H^{'''}_{6}$.  Now assume that $xu\notin A(D)$. By the symmetry between $v$ and $u$, if $xv\in A(D)$, then again we obtain that either $D$ contains a cycle of length 4 or $D\cong H^{'''}_{6}$. We may therefore assume that $A(x\rightarrow  \{u,v\})=\emptyset$. Thus we have $d^+(x,\{u,v,z\})=d^-(x,\{y,z,w\})=0$, i.e., $d(x)\leq 4$, which is a contradiction.\\

\noindent\textbf {Case (b).} There exists a vertex $v\in I(x)$ such that $A(E(x,y)\rightarrow v)= \emptyset$. 

Then
 we have,  $A(E(x,y)\cup O(y)\setminus \{v\}\rightarrow v)=\emptyset$.  By Lemma 3.4(ii), $|E(x,y)\cup O(y)\setminus \{v\}|\leq n$. This together with $|O(y)|\geq n$ implies that $|E(x,y)|=n-2=1$, $n=3$,  $xv\in A(D)$ and $od(v)\geq 3$. If $vy\in A(D)$, then $C_4=xvyv_1x$, where $v_1\in I(x)\setminus \{v\}$. We may assume that $vy\notin A(D)$. Then $v\rightarrow \{z, w\}$, where $E(x,y)=\{w\}$.
 We have that $A(w \rightarrow \{x,y,v\})=\emptyset$. Therefore, $wv_1\in A(D)$, where $\{v,v_1\}=I(x)$, and $C_4=xvwv_1x$, a contradiction.  This completes the proof of Claim 4.1, i.e., 
 $I(x)=O(y)=S(x,y)$.\end{proof}

 Now we divide Case 1 into two subcases.
 
\noindent\textbf {Subcase 1.1.} $ A(x\rightarrow S(x,y))\not= \emptyset$.

 By Claim 4.1, $I(x)=O(y)=S(x,y)$. Let $xu$ be an arbitrary arc in $A(x\rightarrow S(x,y))$. If $uy\in A(D)$, then  $C_4=xuyu_1x$, where $u_1\in S(x,y)\setminus \{u\}$,  a contradiction.
  We may assume that $uy \notin A(D)$. Then, $A(u\rightarrow  S(x,y)\cup \{y\}\setminus \{u\})= \emptyset$.
 This together with $|S(x,y)|\geq n$ and Lemma 3.4(ii) implies that $u\rightarrow E(x,y)$. 
 Note that $|S(x,y)|=od(y)=n$ and $|E(x,y)|= p-n-2$. If for some $w\in E(x,y)$ and $v\in S(x,y)\setminus \{u\}$, $wv\in A(D)$, then $C_4=xuwvx$, a contradiction. We may therefore assume that $A(E(x,y)\rightarrow \{x\}\cup S(x,y)\setminus \{u\})= \emptyset$.
 Then by Lemma 3.4(ii), $E(x,y)\rightarrow y$ and $D[E(x,y)]$ 
 is a complete digraph, which in turn implies that  $|E(x,y)|=1$ since for otherwise, we have $C_4=yuww_1y$, where $w,w_1\in E(x,y)$.
 So, we have $E(x,y)=\{w\}$, $O(w)=\{u,y\}$, i.e.,  $od(w)=2\geq n-1$. This means that $n=2$ or $n=3$. Let $n=2$ and $I(x)=\{u,v\}$. Then it is not difficult to see that
  $id(v,\{u,w\})=0$, $xv\in A(D)$ and 
  $vy\in A(D)$ and $C_4=xvyux$, a contradiction. Let now $n=3$. If 
  $xw\in A(D)$, then $C_4=xwyux$,  a contradiction.  We may therefore assume that $xw\notin A(D)$. From $id(w) =3$ and $A(\{y,x\}\rightarrow w)=\emptyset$ it follows that $S(x,y)\rightarrow w$.    
  Since $D[S(x,y)]$ is an empty digraph and $od(w,S(x,y)\setminus \{u\})=0$, we obtain that $x\rightarrow S(x,y)$.  Therefore, $C_4=xvwux$ is a cycle of length 4 in $D$,   a contradiction. \\

  \noindent\textbf {Subcase 1.2.} $ A(x\rightarrow S(x,y))= \emptyset$.

By the digraph duality, we can assume that $A(S(x,y)\rightarrow y)=\emptyset$.  This together with   $|S(x,y)|\geq n $ and Lemma 3.4(ii) implies that $x\rightarrow E(x,y) \rightarrow y$. Note that $E(x,y)\not= \emptyset$  since $od(x) \geq 2$. Therefore, $C_4=xwyux$, where $w\in E(x,y)$, $u\in S(x,y)$,  a contradiction. This completes the proof of Case 1.\\
 
\noindent\textbf {Case 2.} For each arc $xy\in A(D)$ if $yx\notin A(D)$,  then $od(y)<n$ and $id(x)<n$. 

From the conditions of the  theorem it follows  that $od(y)=id(x)=n-1$,   $p=2n\geq 6$, $id(y)\geq n$ and $od(x)\geq n$. If $S(x,y)= \emptyset$, then using (2) and Lemma 3.4(ii) it is easy to see that
  $D[I(x)\cup \{x\}]$ and $D[O(y)\cup \{y\}]$ are complete digraphs. Therefore, either $n\geq 4$ and $D$ has a cycle of length 4 or $n=3$ and $D\in \cal{H}$$(3,3)$.
  We may therefore assume that 
 $S(x,y)\not= \emptyset$. Since $id(y)|\geq n$ and $od(x)\geq n$, taking into account considered Case 1, we can assume that $O(y)\rightarrow y$ and $x\rightarrow I(x)$.  Hence it is easy to see that $|S(x,y)| =1$, $I(x)\not= O(y)$  and  $|E(x,y)| =1$. Let $E(x,y):=\{w\}$ and $S(x,y):=\{z\}$. From   $A(O(y)\setminus \{z\}\rightarrow I(x)\cup \{x\})=\emptyset$ and 
$A(O(y)\cup \{y\}\rightarrow I(x)\setminus \{z\})=\emptyset$
 it follows that
 $O(y)\setminus \{z \} \rightarrow w \rightarrow I(x)\setminus \{z\}$ and $D[O(y)\cup \{y\}\setminus \{z\}]$ and $D[I(x)\cup \{x\}\setminus \{z\}]$ are complete digraphs. This together with $xz\in A(D)$ and $zy\in A(D)$ implies that $a(z,w)=0$ since $D$ has no cycle of length four. From
 $a(z,w)=0$,
   $d(z)\geq 2n-1\geq 5$ and the fact that
 $A(O(y)\setminus \{z\}\rightarrow z)=A(z\rightarrow I(x)\setminus \{z\})=\emptyset$,  
 it follows that either there exists $v\in I(x)\setminus\{z\}$ or $u\in O(y)\setminus\{z\}$  such that $vz\in A(D)$ or $zu\in A(D)$. Without loss of generality we may assume that $zu\in A(D)$. From this it is easy to see that $A(I(x)\setminus \{z\}\rightarrow  z)=\emptyset$, $wy\notin A(D)$ and $n=3$. Thus we have $O(y)=\{z,u\}$ and $I(x)=\{v,z\}$. From $d^+(w,\{x,y,z\})=d^-(w,\{y,z\})=0$ it follows that $d^+(w)=2$ and $d^-(w)\geq 3$.    Therefore
 the arcs    $wu, vw, xw, vu$ are present in $A(D)$. Therefore, $D$ contains only $Q_6^*=xzyuwvx$  and the arcs $xy$, $xw$, $vu$ and $zu$.  This means that  $D\cong  H^{'}_{6}$. 
We have completed the proof of the theorem. \end{proof}

\section {A sufficient condition for the existence of  a pre-Hamiltonian \\cycle 
in a digraph}

\noindent\textbf{Theorem 5.1.} Let $D$ be a digraph of order $ p\ge5$  with minimum degree  at least $p-1$ and with minimum semi-degree  at least $p/2-1$. Then either $D$ contains a pre-Hamiltonian cycle  or 
$D\in$ \,\,$\cal{H}$$(n,n)$\, \,$\cup \,\,\cal{B}$ $(n,n)$ \,\,$\cup$  $\{[(K_n\cup K_n)+K_1]^*, H_{2n}, H^\prime_{2n}, Q^*_5\}$.

\begin{proof} Suppose, on the contrary, that the theorem is not true. In particular, $D$  contains no cycle of length $p-1$ and $D\notin$ $\cal{H}$$(n,n)$.  Let $ C:= C_{m}:= x_1x_2\ldots x_mx_1$ be an arbitrary non-Hamiltonian cycle of maximum length in $D$. It is not difficult to see that $m\in [3,p-2]$ (if $m=2$, then $D=Q^*_5$).
  From Lemma 3.1 and the maximality of $m$ it follows that for each vertex $y\in B:=V(D)\setminus V(C)$ and for each  $i\in [1,m]$,
$$
 d(y,V(C))\leq m,  \,\,  d(y,B)\geq p-m-1  \,\, \hbox {and if} \,\, x_iy\in A(D), \, \hbox {then } \, yx_{i+1}\notin A(D).
 \eqno{(5)}
$$
 Using the inequality $d(y,B)\geq p-m-1$, it is not difficult to show that the following claim is true. \\

\noindent\textbf {Claim 5.1.}  {\it  If  for  two distinct vertices $u$ and $v$ of $B$  in subdigraph $D[B]$ there is no   $(u,v)$-path, then in $D[B]$  there is a  $(v,u)$-path of length at most 2}. \\

First, we will prove Claims 5.2 and 5.3.\\
  
\noindent\textbf{Claim 5.2.} {\it The induced subdigraph $D[B]$ is strong}.
 
\begin{proof}  Suppose, on the contrary, that $D[B] $ is not strong. Let $D_1$, $D_2$, \ldots  ,    $D_s$ $(s\geq 2)$ be the strong components of $D[B]$,  labeled in such a way that no vertex of $D_i$ dominates a vertex of $D_j$ whenever $i>j$. From Claim 5.1 it follows  that for every pair of vertices $ y\in V(D_1)$ and $z\in V(D_s)$ in  $D[B]$ there is a path  from  $y$ to  $z$ of length 1 or 2. We choose  vertices $y$ from $ V(D_1)$ and $z$ from $V(D_s)$ respectively, such that the path $y_1y_2\ldots y_k$, where $y_1:=y$ and $y_k:=z$, has the minimum length among  all the paths in $D[B]$ with an origin vertex in $D_1$ and a terminal vertex in $D_s$. 
 By Claim 5.1, $k=2$ or $k=3$. We will consider the following three possible cases.\\

\noindent\textbf {Case 5.2.1.} $k<|B|=p-m.$

From the maximality of $C$ it follows that if $x_iy_1\in A(D)$ with $i\in [1,m]$, then $A(y_k\rightarrow \{ x_{i+1},x_{i+2},$  $\ldots ,x_{i+k}\})=\emptyset$. Since $D$ is strong, it follows that $V(C)\not\subseteq I(y_1)$. Therefore the vertex $y_k$ does not dominate at least $id(y_1,V(C))+1$ vertices of $C$. On the other hand, we have $A(y_k\rightarrow V(D_1))=\emptyset$ and $ I(y_1)\subset V(C)\cup V(D_1)$. Hence the vertex $y_k$ does not dominate at least 
$
id(y_1,V(C))+1+ id(y_1,V(D_1))+|\{y_1,y_k\}|= id(y_1)+3
$ 
vertices. From this we obtain $od(y_k)\leq p-id(y_1)-3\leq p/2-2$, which is a contradiction.\\

\noindent\textbf {Case 5.2.2.} $k=|B|=2$. 

It is easy to see that $s=2$, $m=p-2$, $V(D_1)=\{ y_1 \}$, $V(D_2)=\{ y_2 \}$, $ I(y_1)\subset V(C)$ and 
$a^+(x_i, y_1) +  a^+(y_2, x_{i+2}) \leq 1$      
for all $i\in [1,m]$. Hence the vertex $y_2$ does not dominate at least $id(y_1)+|\{y_1,y_2\}|=id(y_1)+2$ vertices. Therefore $od(y_2)\leq p-id(y_1)-2\leq p/2-1$. This together with 
$id(y_1)\geq p/2-1$ and $od(y_2)\geq p/2-1$  implies that $id(y_1)=od(y_2)=p/2-1$. Therefore, $p=2n\geq 6$, $m=p-2\geq4$, $id(y_1)=od(y_2)=n-1$  and
$$
y_2x_i\in A(D) \quad \hbox {if and only if} \quad  x_{i-2}y_1 \notin A(D).      \eqno {(6)} 
$$
By Lemma 3.1, it is easy to see that $d(y_1)=d(y_2)=2n-1$ and $od(y_1)=id(y_2)=n$. Now using Lemma 3.2, we obtain, if $a(y_j,x_i)=0$, then $x_{i-1}y_j$ and $y_jx_{i+1}\in A(D)$. 
We divide this case into two subcases.\\

\noindent\textbf {Subcase 5.2.2.1.} There exists an $i\in [1,m]$ such that $y_1\rightarrow \{ x_i,x_{i+1} \}$. 

 Without loss of generality,  we may assume that  $y_1\rightarrow \{x_2,x_3\}$  and  $a(x_1,y_1)=0$. Therefore, $x_my_1\in A(D)$. From this, (5) and (6) it follows that $x_2y_1 \notin A(D)$, $y_2\rightarrow \{x_3, x_4\}$ and $a(x_2,y_2)=0$. 
  Since $d(y_2,V(C))=2n-2$ and the vertex $y_2$ cannot be inserted into the path $x_3x_4 \ldots x_mx_1$,  by Lemma 3.2, we have $x_1y_2\in A(D)$. If $x_2x_1 \in A(D)$, then $ C_{2n-1}=x_my_1x_2x_1y_2x_4\ldots x_m$. This contradicts our supposition that $D$  contains no cycle of length $p-1$. We may therefore assume that $x_2x_1\notin A(D)$. From this and $a(x_2,y_2)=a^+(x_2,y_1)=0$ it follows that $d(x_2,\{x_3,x_4,\ldots ,x_m \})\geq 2n-3$. Therefore by Lemma 3.2, $x_mx_2 \in A(D)$
 since the vertex $x_2$ cannot be inserted into the path $x_3x_4\ldots x_m$. Now it is easy to see that $a^+(x_i,y_1) + a^+(y_2,x_{i+1}) \leq 1$ for all $i\in [2,m-1]$. 
Therefore, $x_3y_1\notin A(D)$ since $y_2x_4\in A(D)$. From this and (6) it follows that $y_2x_5\in A(D)$ and $x_4y_1\notin A(D)$. Continuing in this manner, we obtain that $A(\{x_4,x_5,\ldots, x_{m-1}\}\rightarrow y_1)=\emptyset $. This together with $A(\{x_1,x_2,x_3\}\rightarrow y_1)=\emptyset$ implies that  $A(\{x_1,x_2,\ldots, x_{m-1}\}\rightarrow y_1)=\emptyset$, i.e., $id(y_1)=1$, which is a contradiction.\\

\noindent\textbf {Subcase 5.2.2.2.} For all $i\in [1,m]$, $|A(y_1\rightarrow \{ x_i,x_{i+1})|\leq 1$. 

Since $od(y_1)=n$, we can assume that $O(y_1)=\{x_1,x_3,\ldots ,x_{2n-3},y_2 \}$. Using this and the fact that $od(y_2)=id(y_1)=n-1$, we obtain $I(y_1)=\{x_1,x_3 ,\ldots , x_{2n-3}\}$. Then by (6),
 $$
O(y_2)=\{x_2,x_4,\ldots ,x_{2n-2}\} \quad \hbox {and} \quad  I(y_2)=\{y_1,x_2,x_4,\ldots ,x_{2n-2}\}.
$$
If for two distinct vertices $x_i,x_j\in \{x_1,x_3,\ldots ,x_{2n-3}\}$,       $x_ix_j\in A(D)$,  then  $C_{2n-1}=y_1x_ix_jx_{j+1}\ldots\\ x_{i-1}y_2x_{i+1}\ldots $ $x_{j-2}y_1$, when  $| \{x_{i+1},x_{i+2}, \ldots ,x_{j-1}\}|  \geq 2$  and $C_{2n-1}=x_ix_jy_1y_2x_{j+1} x_{j+2}\ldots x_{i-1}x_i$,  when $|\{x_{i+1},x_{i+2},$ $\ldots ,$ $x_{j-1}\}|=1$, respectively. This contradicts the supposition that  $C_{p-1}\not\subset D$. We may therefore assume that
 $D[\{x_1,x_3,\ldots ,x_{2n-3},y_2 \}]$ is an empty digraph.
Considering the converse digraph of $D$, by the same arguments we obtain  that
$D[\{x_2,x_4,\ldots ,x_{2n-2},y_1 \}]$  also is an empty digraph.
Therefore  $D\in$ $\cal{B}$${(n,n)}$, which contradicts our supposition that the theorem is not true.\\

 \noindent\textbf {Case 5.2.3.} $k=\mid B\mid =3$.  

From the minimality of  $k$ it follows that $y_1y_3\notin A(D)$, $s=3,$  $A(\{y_2,y_3\} \rightarrow  y_1)=\emptyset $ and $V(D_1)=\{ y_1\}$. Hence $I(y_1)\subset V(C)$. On the other hand, from the maximality of the cycle $C$ it follows that for each $i\in [1,m]$,  
if $x_iy_1\in  A(D)$, then $A(y_2\rightarrow \{x_{i+1},x_{i+2}\})=\emptyset$.
Therefore, $y_2$ does not dominate at least $id(y_1)+1+|\{y_1,y_2\}|=id(y_1)+3$ vertices, which means that $od(y_2)\leq p-id(y_1)-3\leq p/2-2$, a contradiction. This completes the proof that the induced subdigraph $D[B]$ is strong.\end{proof} 

\noindent\textbf {Claim 5.3.} {\em  There exist at least two distinct vertices of $V(C)$  which are adjacent to some vertices of $B$}.

\begin{proof} Assume that Claim 5.3 is not true. Then exactly one vertex of  $C$, say $x$, is adjacent to some vertices of $B$. Therefore, for each vertex $x_i \in V(C)\setminus \{x\}$ and for each vertex $y\in B$ we have 
 $$
d(x_i)=d(x_i,V(C))\leq 2m-2 \quad \hbox {and} \quad d(y)=d(y,B)+d(y,\{x\})\leq 2p-2m.
$$
 Since $d(x_i)+d(y)\ge 2p-2$, we conclude that the inequalities above are equalities. This implies that the subdigraphs $D[V(C)]$ and $D[B\cup \{x\}]$ are complete digraphs. From $d(x_i)=2m-2\geq p-1$ and $d(y)=2p-2m\geq p-1$, we obtain  $p=2m-1$, i.e., $p-m=m-1$. Therefore, $D\cong [(K_{m-1}\cup K_{m-1})+K_1]^*$, which contradicts our supposition. This  proves Claim 5.3. \end{proof}

Since $D$ is strong, we see that  $A(V(C)\rightarrow B)\not=\emptyset$ and $A(B\rightarrow V(C))\not=\emptyset$. This together with Claim 5.3  implies  that there are vertices $x_a\not= x_b$, $x_a,x_b\in V(C)$  and  $x,y\in B$ such that $D$ contains the arcs $x_ax$, $yx_b$. Without loss of generality, we can assume that   $|V(C[x_a,x_b])|$ is as small as possible.  Note that       if $x_b\not=x_{a+1}$, then
$$
A(\{x_{a+1},x_{a+2},\ldots ,x_{b-1}\},B)=\emptyset.     \eqno {(7)} 
$$
For the sake of clarity, assume that $x_b:=x_1$ and $x_a:=x_{m-h}$, where  $0\leq h\leq m-2$. We  consider the following two cases, depending  on whether $x_1=x_{m-h+1}$ or not.\\

\noindent\textbf {Case 1.} $x_{m-h+1}\not=x_1$, i.e., $h\geq 1$.

Consider the paths $P_0$, $P_1$, \ldots , $P_k$ (where $0\leq k\leq h$ and $k$ is as large as possible), where $P_0:=P= x_1x_2\ldots x_{m-h}$ and  each  $P_i$ path, where $i\in [1,k]$, is an extended path obtained from the path $P_{i-1}$ with a vertex $z_i\in \{x_{m-h+1}, x_{m-h+2},\ldots , x_m \}\setminus \{z_1,z_2, \ldots, z_{i-1} \}$. 
Note that the path $P_i$, $i\in [0,k]$, contains $m-h+i$ vertices. We may assume that the extended  path $P_k$ does not contain some vertices $y_1,y_2,\ldots ,y_d \in \{x_{m-h+1},x_{m-h+2},\ldots ,x_{m}\}$, where $1\leq d\leq h$. Therefore, using (7) and Lemma 3.2, for each $z\in B$ and for each $y_i$ we obtain
$$
d(z)=d(z,B)+d(z,V(C))\leq 2p-2m-2+m-h+1=2p-m-h-1
$$
and
$$
d(y_i)=d(y_i,V(C))\leq m-d+1+2d-2= m+d-1.
$$
Hence it is clear that $$2p-2\leq d(z)+d(y_i)\leq 2p+d-h-2.$$ 
  It is not difficult to see that $h=d$, $d(z,V(C))=m-h+1$, $d(y_i,V(C))=m+h-1 $ and both the subdigraphs  $D[B]$ and $D[\{x_{m-h+1},x_{m-h+2},\ldots ,x_{m}\}]$ are  complete digraphs.
  By Lemma 3.2(ii), we also have $x_{m-h} \rightarrow B \cup \{x_{m-h+1},x_{m-h+2},\ldots ,x_{m}\} \rightarrow x_1 $. Using the inequalities $p-1\leq d(y_i)\leq m+h-1$ and $p-1\leq d(z)\leq 2p-m-h-1$, it is easy to see that $h=|B|=p-m \geq 2$ and the path $P=x_1x_2\ldots x_{m-h}$ has the maximum length among all the paths from $x_1$ to $x_{m-h}$ in the subdigraphs $D[V(C)]$ and  $D[B \cup \{x_1 ,x_2  ,\ldots ,x_{m-h}\}]$. Therefore by Lemma 3.3, there exist integers $l\in [1,m-h]$ and $ r\in [1,m-h] $ such that for all $u\in B$ and for all $z\in \{x_{m-h+1},x_{m-h+2},\ldots , x_{m}\}$ the
 following holds 
$$
\left\{\begin{array}{rcl}
O(u,V(P))=\{x_1,x_2,\ldots ,x_l\}, \quad I(u,V(P))=\{x_l,x_{l+1},\ldots ,x_{m-h}\}, \\
O(z,V(P))=\{x_1,x_2,\ldots ,x_r\}, \quad I(z,V(P))=\{x_r,x_{r+1},\ldots ,x_{m-h}\}.  \\ 
\end{array}
\right.
\eqno (8)
$$
According to the digraph duality, we can assume that   $l\leq r$. (See Figure 1 (b) for an illustration).

 Let $l=1$ and let $u$ be an arbitrary vertex of $B$. Then from $od(u)\geq p/2-1$ and (8) it follows that $h\geq p/2-1$ and $ p\geq  2(p/2-1)+m-h=p-2+m-h$. So, $m-h=2$ since $m-h\geq  2$. We see that $h=p/2-1$, $p=2n\geq 6$,  $h=n-1$, and $r=2$ since $id(x_2)\geq n-1\geq 2$. Therefore $D\in \{H_{2n},H'_{2n}\}$, which contradicts  our  supposition. 

 Let now $l\geq 2$. By the digraph duality, we can assume that $r\leq  m-h-1$.    Since $D[B]$ and $D[\{x_{m-h+1},x_{m-h+2},\ldots ,x_{m}\}]$  are complete digraphs and (8), for each vertex  $ z\in \{x_{m-h+1},x_{m-h+2},\\\ldots ,x_m \}$  we have $I(z)=\{x_r,x_{r+1},\ldots ,x_m \}\setminus \{ z\}$. This implies that $m-r\geq p/2-1$. If $i\in [r+1,m-h]$  and $x_1x_i\in A(D)$ then by (8) and $2\leq l\leq r\leq m-h-1$, we have $ C_{m+1}=x_1x_ix_{i+1}\ldots x_mx_2\ldots x_{i-1}xx_1$, where $x\in B$, a contradiction. Because of this and $2\leq l\leq r$, we may assume that
$
A(x_1\rightarrow B\cup \{ x_{r+1},x_{r+2},\ldots , x_m \})=\emptyset.
$
Therefore, since $m-r\geq p/2-1$ and $|B|=h\geq 2$, we obtain $od(x_1)\leq  p-1-h-(m-r)\leq p/2-h\leq p/2-2$, which contradicts  that $od(x_1)\geq p/2-1$. This contradiction completes the discussion of Case 1.\\

\noindent\textbf {Case 2.} $x_{m-h+1}=x_1$, i.e., $h=0$.

Recall that $x,y \in B$, $x_{m-h}x=x_mx\in A(D)$ and $yx_1\in A(D)$.      Therefore, any path from  $x$ to  $y$ in   $D[B]$ is a Hamiltonian path.   Let $ u_1u_2\ldots u_{p-m}$ be an $(x,y)$-Hamiltonian path in $D[B]$, where $u_1:=x$, $u_{p-m}:=y$. Observe that if $1\leq i<j\leq p-m$, then  $u_iu_j\in A(D)$ if and only if $j=i+1$.
For  Case 2 ($h=0$), we will first prove a number of claims (Claims 5.4-5.8).\\

\noindent\textbf {Claim 5.4.} $p-m=2$, i.e., $m=p-2$.

\begin{proof} Suppose, to the contrary, that $p-m\geq 3$.  From the above observation it follows  that $u_1u_{p-m}\notin A(D)$ and $od(u_1,B)=id(u_{p-m},B)=1$. Using this and (5), we obtain
$$
p-1\leq d(u_1)\leq m+1+id(u_1,B) \,\,  \hbox {and} \,\, p-1\leq d(u_{p-m})\leq m+1+od(u_{p-m},B).
$$
Therefore, $id(u_1,B)\geq p-m-2$  and  $od(u_{p-m},B)\geq p-m-2$. Therefore, in  $D[B]$ there is a path from  $u_{p-m}$ to  $u_1$ of length $k=1$ or $k=2$, for otherwise $p-m\geq 2p-2m-2$, which contradicts that  $p-m\geq 3$. For any integer $l\geq 1$, put
$$
I^+_l(u_{p-m}):=\{x_j |\,x_{j-l}u_{p-m}\in A(D) \}.
$$
Since  $id(u_{p-m},V(C))=id(u_{p-m})-1$ and \, $V(C)\not\subseteq I(u_{p-m})$, we see that for each $l\in [1,2]$, 
$$
| I^+_l(u_{p-m})\cup I^+_{l+1}(u_{p-m}) | \geq id(u_ {p-m}).
$$
From the maximality of the cycle $C$ it follows that
$
A(u_1\rightarrow I^+_k(u_{p-m})\cup I^+_{k+1}(u_{p-m}))=\emptyset.
$
This together with  $A(u_1\rightarrow \{u_3,u_4,\ldots ,$ $u_{p-m}\})=\emptyset$ implies that 
$$
p/2-1\leq od(u_1)\leq p-1-| I^+_k(u_{p-m}) \cup I^+_{k+1}(u_{p-m})| -(p-m-2)$$ $$\leq m+1-id(u_{p-m})\leq m+1-p/2+1=m+2-p/2.
$$
Therefore, since $m\leq p-3$, we obtain  that $p=2n$,  $p-m=3$  and $od(u_1)=id(u_3)=n-1$. Hence,  $id(u_1)\geq n$  and $od(u_3)\geq n$. 
We now claim that $D$ contains the arcs $u_3u_1$ and $u_2u_1$.  Indeed, for otherwise, $id(u_1,V(C))\geq n-1$ and if $x_iu_1\in A(D)$, then $od(u_3, \{x_{i+2},x_{i+3}\})=0$. From this it is not difficult to see that  $od(u_3)\leq n-1$, which contradicts the fact that   $od(u_3)\geq n$.
Similarly, we can see that $u_3u_2\in A(D)$. So, the arcs $u_3u_1$, $u_2u_1$, $u_3u_2$ are in $A(D)$. Since $id(u_3,V(C))=n-2$, $m\geq n$,  $m\geq id(u_3,V(C))+2$ and $C$ is a  non-Hamiltonian cycle of maximal length, it follows that $|\cup _{i=1}^3 I^+_i(u_{p-m})| \geq n$ and $od(u_1, \cup _{i=1}^3 I^+_i(u_{p-m}))=0$. This together with $u_1u_3\notin A(D)$  implies that  $od(u_1)\leq n-2$, a contradiction. This completes the proof that $p-m=2$. \end{proof}

Let $B:=\{u,v\}$. Then by the Claim 5.2 we have  $uv\in A(D)$ and $vu \in A(D)$.\\ 

\noindent\textbf {Remark 2.} Claims 5.5–5.8 also hold for vertex $v$ by virtue of the symmetry between $u$ and $v$.\\

\noindent\textbf {Claim 5.5.} If $x_iu$ and $ux_{i+2} \in A(D)$ with $i\in [1,m]$, then $a(x_{i+1},v) =2$, i.e., $x_{i+1}v$ and $vx_{i+1}\in A(D)$.

\begin{proof} Notice that the cycle $x_iux_{i+2}x_{i+3}\ldots x_i$ has length $m$ and the vertices $v$ and $x_{i+1}$ are not on this cycle.  By Claim 5.2,  the subdigraph $D[\{v,x_{i+1}\}]$ is strong. Therefore $vx_{i+1}$ and $x_{i+1}v \in A(D)$. \end{proof}

Using Claim 5.5, the facts that $uv$, $vu\in A(D)$ and $C$  is a   non-Hamiltonian cycle of maximum length in $D$, we obtain the following claim.

\noindent\textbf {Claim 5.6.} If $i\in [1,m]$, then
$$
| A(\{x_i,x_{i+1}\}\rightarrow u)| +| A(u\rightarrow x_{i+3})| \leq 2 \,\, \hbox {and}\,\, | A(x_{i-2}\rightarrow u)| + | A(u\rightarrow \{x_i,x_{i+1} \})| \leq 2. 
$$

\noindent\textbf {Claim 5.7.} If $k\in [1,m]$, then  
$|A(\{x_{k-1},x_k\}\rightarrow u)| \leq 1$ and $|A(u\rightarrow \{x_{k-1},x_{k}\})| \leq 1$.

\begin{proof} By the digraph duality, it is enough to prove that if $k\in [1,m]$, then $|A(\{x_{k-1},x_k\}\rightarrow u)| \leq 1$.
Suppose, to the contrary, that there exists  $k\in [1,m]$ such that $\{ x_{k-1},x_k \}\rightarrow u$. Without loss of generality, we may assume that $a(u,x_{k+1})=0$. 
For the sake of clarity, assume that  $x_{k+2}=x_1$. Then  $x_m=x_{k+1}$, $x_{m-1}=x_k$, $x_{k-1}=x_{m-2}$, and by Claim 5.6,  $ux_1\notin A(D)$. 

  Let us first assume that $x_1u\in A(D)$. It is easy to see that  $p\geq 6$, $m\geq 4$ and  
  $od(u, \{x_{m-1},x_m,\\x_1,x_2\})=0$. Using this and $od(u)\geq p/2-1\geq 2$, we see that  $ A(u\rightarrow \{x_3,x_4,\ldots , x_{m-2}\})\not =\emptyset$, which in turn implies that $m-2-2\geq 1$, i.e., $m\geq 5$.   If   $m=5$, then $p=7$ and $od(u)\geq 3$. We know that $od(u, \{x_{m-1},x_m,x_1,x_2\})=0$ and hence, $od(u)\leq 2$, a contradiction. We may therefore assume that $m\geq 6$.  Now, we need to show that, for every  $j\in [3,m-3]$,
$$
 od(u, \{x_j,x_{j+1}\})\leq 1. \eqno {(9)}
$$

Assume that (9) is not true. Then there exists $j\in[3,m-3]$ such that $u\rightarrow \{ x_j,x_{j+1}\}$.
Assume that $j$ is as small as possible. Then  $a(u,x_{j-1})=0$,  and  $x_{j-2}u \notin A(D)$ by Claim 5.6. Hence $j\geq 4$. Since the vertex $u$ cannot be inserted into the cycle $C$, $ux_1 \notin A(D)$  and  $x_{j-2}u \notin A(D)$, using Lemma 3.2, we obtain
$
p-1\leq d(u)=d(u,\{x_1,x_2,\ldots ,x_{j-2}\})+ d(u,\{x_j,x_{j+1},\ldots , x_{m-1}\}) +d(u,\{v\})$$ $$\leq j-3+m-1-j+4=m=p-2,
$
which is a contradiction. This proves (9).\\
 
Recall that $A(u\rightarrow  \{x_{m-1},x_m,x_1,x_2\})=\emptyset$. From this and (9) it follows that: If $|\{x_3,x_4,\ldots ,\\ x_{m-2}\}|=p-6$ is even, then $p/2-1\leq od(u)\leq (p-6)/2+1$\, a contradiction.  If $|\{x_3,x_4,\ldots ,\\ x_{m-2}\}|=p-6$ is odd, then $p$ is odd and $(p-1)/2\leq od(u)\leq (p-7)/2+2$, a contradiction. 

Suppose next that $x_1u\notin A(D)$. Then by Claim 5.6, $a(u,x_1)=0$.  Notice that $m\geq 4$ and $d(u,\{x_2,x_3, \ldots,x_{m-1}\})$ $\geq p-3$. Hence by Lemma 3.2(ii),  $ ux_2\in A(D)$. We have $a(v,x_m)=0$ and $vx_1\notin A(D)$. By Lemma 3.2(iii), it is easy to see that  $x_{m-1}v\in A(D)$ and  $d(v,\{x_1,x_2,\ldots ,x_{m-1}\})=p-3$. If $x_1v\notin A(D)$, then $a(v,x_1)=0$, and by Lemma 3.2, $vx_2\in A(D)$. Now we have $x_{m-1}u$, $vx_2\in A(D)$ and $A(\{u,v\},\{x_m,x_1\})=\emptyset$. Thus we have the considered Case 1 ($h\geq 1$). We may therefore assume that $x_1v\in A(D)$. We also may assume that $x_{m-2}v\notin A(D)$  (for otherwise $\{x_{m-2},x_{m-1},x_1\}\rightarrow v$ and for the vertex $v$ the considered case $x_1u\in A(D)$ holds). Since $v$ cannot be inserted into the path $x_1x_2\ldots x_{m-2}$, using the facts that $x_{m-2}v\notin A(D)$,  $vx_1\notin A(D)$  and Lemma 3.2(iii), we obtain that $d(v,\{x_1,x_2,\ldots , x_{m-2}\})\leq m-3=p-5$. Then from $a(x_m,v)=0$ and $d(v)\geq p-1$ it follows that $vx_{m-1}\in A(D)$. If
$x_mx_2\in A(D)$, then $C_{m+1}=x_mx_2x_3\ldots x_{m-2}uvx_{m-1}x_m$, a contradiction. If $x_{m-2}x_m\in A(D)$, then $C_{m+1}=x_{m-2}x_mx_1vux_2x_3\ldots x_{m-2}$, a contradiction. We may therefore assume that $x_mx_2\notin A(D)$ and $x_{m-2}x_m\notin A(D)$.
Note that $x_m$ cannot be inserted into the path $x_2x_3\ldots x_{m-1}$ (for otherwise there exists an $(x_2,x_{m-1})$-path, say $R$, with vertex set $\{x_2,x_3,\ldots , x_m\}$, and hence $x_{m-1}vuRx_{m-1}$ is a cycle of length $p-1$, a contradiction). From this, $d(x_m,\{u,v\})=0$, $d(x_m)\geq p-1$ and Lemma 3.2(iii) it follows that  $p-1\le d(x_m)= d(x_m,\{u,v\})+d(x_m,\{x_{m-1},x_1\})+d(x_m,\{x_2,x_3,\ldots , x_{m-2}\})\le 4+m-4=m=p-2$,
 a contradiction. Claim 5.7 is proved. \end{proof}

\noindent\textbf {Claim 5.8.} If  $k\in [1,m]$, then $a^+(x_k,  u) +a^+(u, x_{k-1}) \leq 1$.

\begin{proof} Suppose, to the contrary, that there exists  $k\in [1,m]$ such that  $x_ku$, $ux_{k-1}\in A(D)$. For the sake of clarity, let us assume that  $x_k:=x_2$. From Claim 5.7  and (5) it follows that
$d(u,\{x_m,x_3\})=a^+(x_1, u)=a^+(u, x_2)=0$. First we show that $m\geq 5$. Assume that this is not the case. Then either $m=4$ or $m=3$. Let $m=4$, i.e., $p=6$. Using Claim 5.7 and (5), we obtain that $d(u,\{x_3,x_4\})=0$, which together with $a^+(x_1, u)=a^+(u, x_2)=0$ implies that $d(u)=4$, a contradiction. Let now $m=3$, i.e., $p=5$. From $a^+(x_1, u)=0$ it follows that $D$ is not a symmetric digraph. Therefore, $D\notin \{Q^*_5, [K_2\cup K_2]+K_1]^*\}$. Then by Theorem 4.3(ii), $D$ contains a cycle of length $4=p-1$, which contradict our initial  supposition. So, $m\geq 5$.
It is easy to see that  $d(u,\{x_4,x_5,\ldots ,x_{m-1}\})\geq p-5$. Then, since the vertex $u$ cannot be inserted into the path $x_4x_5\ldots x_{m-1}$, from Lemma 3.2 it follows that  $ux_4\in A(D)$ and $x_{m-1}u\in A(D)$. Then by Claim 5.5,  $a(v,x_3) =a(v,x_m) =2$.   Therefore by Claim 5.7, $d(v,\{x_1,x_2\})=0$. Since  $D$ contains no  cycle of length $m+1$, it is not difficult to see that $x_mx_2$, $x_2x_4$ and $x_3x_2$ are not in $A(D)$ (if $x_3x_2\in A(D)$, then $C_{m+1}=x_mvx_3x_2ux_4x_5\ldots x_m$). So we have $d(x_2,\{u,v,x_1,x_3\})\leq 4$ and $d(x_2,\{x_4,x_5,\ldots ,x_m\})\geq p-5$. Therefore, since neither $x_mx_2$ nor $x_2x_4$ are in $A(D)$, using Lemma 3.2(iii), we can insert the vertex $x_2$ into the path $x_4x_5\ldots x_m$, 
 i.e., for some $i\in [4,m-1]$, $x_ix_2$,  $x_2x_{i+1} \in A(D)$, and obtain a cycle $x_mvx_3x_4\ldots x_ix_2x_{i+1}\ldots x_m$  of length $m$, which does not contain the vertices $x_1$ and $u$. By Claim 5.2, the subdigraph $D[\{x_1,u \}]$ must be strong. But this is not the case since    $x_1u\notin A(D)$. This contradiction completes the proof of Claim 5.8. \end{proof}

We now divide Case 2 ($h=0$) into two subcases depending on whether $p$ is odd or not.\\

\noindent\textbf {Subcase 2.1.} \, $ p=2n+1$.

Then $m=2n-1$, $od(u,V(C))\geq n-1$ and $id(u,V(C))\geq n-1$.  From  (5), and Claims 5.7 and 5.8 it follows that each of  vertices $u$ and $v$ is adjacent to at most one of the two consecutive vertices of the cycle $C$. Then it is not difficult to see that $od(u,V(C))=od(v,V(C))=n-1$, $ O(u,V(C))=I(u,V(C))$ and $O(v,V(C))=I(v,V(C))$. Therefore, without loss of generality, we may assume that 
$$
d(u,\{x_2,x_3\})=0 \quad \hbox {and} \quad O(u)=I(u)=\{x_1,x_4,x_6,\ldots , x_{p-3}=x_{2n-2},v\}.  
$$

Let  $m=3$, i.e., $p=5$. If  $vx_3\in A(D)$, then $C_4=x_3x_1uvx_3$, a contradiction. Assume that $vx_3\notin A(D)$ and $x_3v\in A(D)$. Then by Claims 5.7 and 5.8, $d(v,\{x_1,x_2\})=0$. Therefore, $d(v)\leq 3$, which contradicts that $d(v)\geq 4$. We may therefore assume that $a(v,x_3)=0$. Then by Claim 5.3, we have that $a(v,x_2)\not=0$. If $x_2v\in A(D)$, then $C_4=x_1x_2vux_1$, a contradiction. Therefore, we may assume that $x_2v\notin A(D)$. Then $vx_2\in A(D)$ and by Claim 5.7, $a(v,x_1)=0$. Therefore, $d(v)=a(v,u)+a(v,x_2)\leq 3$,  
 a contradiction.
 
Let now $m=2n-1\geq 5$. Since $x_{m-1}u\in A(D)$, by Claim 5.5 we have  $ a(x_m,v) =2$. Then by Claims 5.7 and 5.8, $d(v,\{x_{m-1},x_1\})=0$. Since 
$a(u,x_1)=a(u,x_4) =2$, it follows that neither $vx_3$ nor $x_2v$ are in $A(D)$. This together with $O(v,V(C))=I(v,V(C))$  implies that 
$d(v,\{x_2,x_3\})=0$. Thus, we have $d(v,\{x_1,x_2,x_3\})=0$, which is a contradiction.
 This contradiction completes the discussion of case $p=2n+1$.\\ 

\noindent\textbf {Subcase 2.2.}  $ p=2n $.

  From $d(u)\geq 2n-1$ it follows that either $od(u)\geq n$ or $id(u)\geq n$. By the digraph duality, assume that $od(u)\geq n$.  Now using Claims 5.7 and 5.8, we may assume that 

$$
\left\{\begin{array}{rcl}
u\rightarrow \{x_1,x_3,\ldots ,x_{2n-3}\}, \quad A(u,\{x_2,x_4,\ldots ,x_{2n-2}\})=\emptyset,  \\
I(u)\subseteq \{v,x_1,x_3,\ldots ,x_{2n-3}\}. \\ 
\end{array}
\right.
\eqno (10)
$$
Since $id(u)\geq n-1$, without loss of generality,  we may assume that $\{x_1,x_3,\ldots ,x_{2n-5}\}\rightarrow u$. From this, (10) and Claim 5.5 it follows that for each $i\in [1,n-2]$,
$$
 a(u,x_{2i-1}) =a(v,x_{2i})=2.      \eqno {(11)}
$$
Then by Claims 5.7 and 5.8, we have $d(v,\{x_1,x_3,\ldots ,x_{2n-3}\})=0$. Therefore  $a(v,x_{2n-2})\not=0$ since $d(v)\geq 2n-1$, i.e., either  $vx_{2n-2}\in A(D)$ or $v_{2n-2}v\in A(D)$. If $vx_{2n-2}\in A(D)$, then using Claim 5.5, we obtain $x_{2n-3}u\in A(D)$ and $x_{2n-2}v\in A(D)$.  So, in any case we can assume that $x_{2n-2}v\in A(D)$. Then $id(v)\geq n$ and $\{x_2,x_4,\ldots , x_{2n-2}\}\rightarrow v$.\\

Let us now prove that $D[\{ x_1,x_3,\ldots ,x_{2n-3}\}]$ is an empty digraph.

 Assume that  this is not the case, i.e., 
 for some distinct vertices $x_i, x_j \in \{ x_1,x_3,\ldots ,$ $x_{2n-3}\}$, $x_ix_j\in A(D)$. If  $|\{x_{i+1}, x_{i+2},\ldots ,x_{j-1}\}| =1$, then from (11) and  $x_{2n-2}v\in A(D)$ it follows that:  If $j=2n-3$, then  $i=2n-5$ and   $C_{m+1}=x_{2n-5}x_{2n-3}x_{2n-2}vux_1x_2\ldots x_{2n-5}$.  If $j\not= 2n-3$, then  $C_{m+1}=x_ix_juvx_{j+1}\ldots x_{i-1}x_i$. We may therefore assume that  $|\{x_{i+1},x_{i+2},\ldots ,x_{j-1}\}| \geq 2$. Then $m\geq 6$. Using the above arguments, we obtain: If $i=2n-3$, then $C_{m+1}=x_{2n-3}x_j\ldots x_{2n-5}ux_1\ldots x_{j-1}vx_{2n-4}\\x_{2n-3}$. 
 If $i\not= 2n-3$ and $j\not= 1$, then the arcs $x_{i-1}v$, $vx_{i+1}$,  $x_{j-2}u$  are in $A(D)$, and hence
$C_{m+1}=x_ix_jx_{j+1}\ldots x_{i-1}vx_{i+1}\ldots x_{j-2}ux_i$. Let now $j=1$, then $C_{m+1}=x_ix_1ux_{i+2}\ldots x_{2n-2}vx_{2}x_3\ldots x_i$.  So, in each case we have that $D$ contains a cycle of length  $m+1$, which is a contradiction.  This proves that $D[\{ x_1,x_3,\ldots ,x_{2n-3}\}]$ is an empty digraph.\\

By a similar argument, one can show that $D[\{ x_2,x_4,\ldots ,x_{2n-2}\}]$ is  also  an empty digraph.
Therefore,
$$
A(B[\{v, x_1,x_3,\ldots , x_{2n-3}\}])=A(B[\{u, x_2,x_4,\ldots , x_{2n-2} \}] )=\emptyset.$$

Then by Remark 1,  $D\in$\,$\cal{B}$${(n,n)}$,
which contradicts our supposition. The discussion of Case 2 is completed  and Theorem 5.1 is proved. \end{proof}

It is worth mentioning that, for several sufficient conditions for a digraph to be Hamiltonian, it has been shown that, under these conditions, the digraph also contains a pre-Hamiltonian cycle, except some cases. These can be seen in \cite{[12]},  \cite{[13],[14]}, \cite{[15]}   and the references therein, among others. In addition, a number of papers have considered the existence of the pre-Hamiltonian cycles  in balanced bipartite digraphs. For example, see \cite{[8]} and \cite{[16]}.\\

\textbf{Acknowledgements}

 I would like to thank the anonymous reviewers for their valuable comments and suggestions. I am grateful to Professor Gregory Gutin, who motivated me to present the complete proofs of my earlier unpublished work. Also thanks to Dr. Parandzem Hakobyan for formatting the manuscript of this paper.

\end{document}